\documentclass[11pt]{article}
\usepackage{mathrsfs}
\usepackage{amssymb}
\usepackage{amsmath}
\usepackage[all]{xy}

\def\Ker{\mathop{\rm Ker}\nolimits}
\def\Coker{\mathop{\rm Coker}\nolimits}
\def\Tr{\mathop{\rm Tr}\nolimits}
\def\mod{\mathop{\rm mod}\nolimits}

\title{\Large \bf When Are Torsionless Modules Projective?
}
\author{Rong Luo and Zhaoyong Huang\thanks{\it E-mail: huangzy@nju.edu.cn}\\
{\small \it Department of Mathematics, Nanjing University,
Nanjing 210093, P. R. China}\\
}
\usepackage{amssymb}
\date{}
\begin{document}
\baselineskip=18pt \maketitle

\begin{abstract} In this paper, we study the problem when a finitely generated
torsionless module is projective. Let $\Lambda$ be an Artinian local
algebra with radical square zero. Then a finitely generated
torsionless $\Lambda$-module $M$ is projective if ${\rm
Ext^1_\Lambda}(M,M)=0$. For a commutative Artinian ring $\Lambda$, a
finitely generated torsionless $\Lambda$-module $M$ is projective if
the following conditions are satisfied: (1) ${\rm
Ext}^i_{\Lambda}(M,\Lambda)=0$ for $i=1,2,3$; and (2) ${\rm
Ext}^i_{\Lambda}(M,M)=0$ for $i=1,2$. As a consequence of this
result, we have that for a commutative Artinian ring $\Lambda$, a
finitely generated Gorenstein projective $\Lambda$-module is
projective if and only if it is selforthogonal.
\end{abstract}

\vspace{0.5cm}

\centerline{\large \bf 1. Introduction}

M. Ramras in [G] p.380 raised an open question: For a left and right
Noetherian ring $\Lambda$, when is every finitely generated
reflexive $\Lambda$-module projective? He proved in [R] that if
$\Lambda$ is a commutative Noetherian local ring and $M$ is a
finitely generated $\Lambda$-module such that the sequence of Betti
numbers of $M$ is strictly increasing, then the condition $M$ is
torsionless with ${\rm Ext}^1_{\Lambda}(M,\Lambda)=0$ implies $M$ is
projective. Menzin in [M] proved that if $\Lambda$ is an Artinian
local algebra with radical square zero, then for $\Lambda$ not
Gorenstein all finitely generated reflexive modules are projective.
Recently, Braun in [B] proved that for a commutative Noetherian ring
$\Lambda$, a finitely generated $\Lambda$-module $M$ is projective
if it satisfies the following conditions: (1) The projective
dimension of $M$ is finite; (2) ${\rm End}_{\Lambda}(M)$ is a
projective $\Lambda$-module; and (3) $M$ is reflexive or ${\rm
Ext}^1_{\Lambda}(M,M)=0$. In this paper, we will study a stronger
problem: When is a finitely generated torsionless module projective?

As a common generalization of the notion of projective modules,
Auslander and Bridger in [AuB] introduced the notion of finitely
generated modules of Gorenstein dimension zero. Such a kind of
modules is called Gorenstein projective following Enochs and Jenda's
terminology in [EJ]. It is well known that a projective module is
Gorenstein projective. Then it is natural to ask when the converse
holds true, or equivalently, what is the difference between the
projectivity and Gorenstein projectivity of modules? In views of the
properties of projective modules and Gorenstein projective modules,
we conjecture that the difference between these two classes of
modules is the selforthogonality of modules.

\vspace{0.2cm}

{\bf Gorenstein Projective Conjecture (GPC)} Over an Artinian
algebra, a finitely generated Gorenstein projective module $M$ is
projective if and only if it is selforthogonal.

\vspace{0.2cm}

It is trivial that the necessity in {\bf GPC} is always true. So the
sufficiency is essential in {\bf GPC}. Observe that {\bf GPC} is
related to the question mentioned above. On the other hand, part of
motivation for studying {\bf GPC} is that it is a special case of
the well-known generalized Nakayama conjecture ({\bf GNC}) (it still
remains open), which states that for an Artinian algebra $\Lambda$
and a finitely generated $\Lambda$-module $M$, the condition ${\rm
Ext}_{\Lambda}^i(M\bigoplus \Lambda , M\bigoplus \Lambda)=0$ for any
$i \geq 1$ implies $M$ is projective (see [AuR1]). In this paper, we
will prove that {\bf GPC} is true if $\Lambda$ is commutative, that
is, if $\Lambda$ is a commutative Artinian ring.

In Section 2, we collect some known facts for later use. In Section
3, we prove that for an Artinian local algebra $\Lambda$ with
radical square zero, a finitely generated torsionless
$\Lambda$-module $M$ is projective if ${\rm Ext^1_\Lambda}(M,M)=0$.
For any Artinian algebra, we also give some criteria for judging an
indecomposable torsionless module being projective. In particular,
we provide some support to {\bf GNC}. In Section 4, we prove that if
$\Lambda$ is a commutative Artinian ring, then a finitely generated
torsionless $\Lambda$-module $M$ is projective provided that the
following conditions are satisfied: (1) ${\rm
Ext}^i_{\Lambda}(M,\Lambda)=0$ for $i=1,2,3$; and (2) ${\rm
Ext}^i_{\Lambda}(M,M)=0$ for $i=1,2$. As an immediate consequence,
we have that for a commutative Artinian ring $\Lambda$, a finitely
generated Gorenstein projective $\Lambda$-module is projective if
and only if it is selforthogonal, that is, {\bf GPC} is true for
commutative Artinian rings.

\vspace{0.5cm}

\centerline{\large \bf 2. Preliminaries}

\vspace{0.2cm}

In this section, we give some notions and notations in our
terminology and collect some facts for later use. For a ring
$\Lambda$, we use $\mod \Lambda$ and $J(\Lambda)$ to denote the
category of finitely generated left $\Lambda$-modules and the
Jacobson radical of $\Lambda$, respectively. We use $(-)^*$ to
denote ${\rm Hom}_{\Lambda}(-, \Lambda)$. All modules considered are
finitely generated.

Let $\Lambda$ be an Artinian algebra and
$$P_1 \buildrel {f} \over \to P_0 \to M \to 0$$ a minimal
projective resolution of a module $M$ in $\mod \Lambda$. We call
$\Coker f^*$ is the {\it transpose} of $M$, and denote it by $\Tr
M$. Let $M\in \mod \Lambda$ and $\sigma _{M}: M \to M^{**}$ defined
by $\sigma _{M}(x)(f)=f(x)$ for any $x\in M$ and $f\in M^{*}$ be the
canonical evaluation homomorphism. $M$ is called {\it torsionless}
if $\sigma _{M}$ is a monomorphism; $M$ is called {\it reflexive} if
$\sigma _{M}$ is an isomorphism (see [AuB]). By [Au, Proposition
6.3], we have an exact sequence:
$$0\to {\rm Ext}_{\Lambda^{op}}^1(\Tr M,\Lambda^{op}) \to
M \buildrel {\sigma _M} \over \longrightarrow M^{**} \to {\rm
Ext}_{\Lambda^{op}}^2(\Tr M,\Lambda^{op}) \to 0.$$ On the other
hand, it is easy to see that $\Tr\Tr M$ and $M$ are projectively
equivalent. So, we have that $M$ (resp. $\Tr M$) is torsionless if
and only if ${\rm Ext}_{\Lambda^{op}}^1(\Tr M,\Lambda^{op})=0$
(resp. ${\rm Ext}_{\Lambda}^1(M,\Lambda)=0$); and $M$ (resp. $\Tr
M$) is reflexive if and only if ${\rm Ext}_{\Lambda^{op}}^i(\Tr
M,\Lambda^{op})=0$ (resp. ${\rm Ext}_{\Lambda}^i(M,\Lambda)=0$) for
$i=1, 2$.

We use $\mod _P\Lambda$ to denote the subcategory of $\mod \Lambda$
consisting of modules without non-zero projective summands. For $M$
and $N$ in $\mod \Lambda$, we use $\underline{{\rm Hom}}_\Lambda(M,
N)$ (resp. $\overline{{\rm Hom}}_\Lambda(M, N)$) to denote the set
of the equivalence classes of module homomorphisms modulo those
factoring through a projective (resp. injective) $\Lambda$-module.
For an Artinian algebra $\Lambda$, we denote by $\mathbb{D}$ the
ordinary duality, that is, $\mathbb{D}(-)={\rm Hom}_R(-, R/J(R))$ of
$\Lambda$, where $R$ is center of $\Lambda$ which is a commutative
Artinian ring.

\vspace{0.2cm}

{\bf Lemma 2.1} ([AuR2, Theorem 3.3]) {\it Let} $\Lambda$ {\it be an
Artinian algebra,} $M\in \mod _P\Lambda$ {\it and} $X\in \mod
\Lambda$. {\it Then there is an isomorphism:}
$$\overline{{\rm Hom}}_\Lambda (X, \mathbb{D}\Tr M)\to {\rm Hom}
_{\underline{{\rm End}}(M)^{op}}({\rm Ext}^1_{\Lambda}(M, X), {\rm
Ext}^1_{\Lambda} (M, \mathbb{D}\Tr M)).$$

\vspace{0.2cm}

Recall from [AF] that a module $M$ in $\mod \Lambda$ is called {\it
faithful} if the annihilator of $M$ in $\Lambda$ is zero.

\vspace{0.2cm}

{\bf Lemma 2.2} ([AF, p.217]) {\it Let} $\Lambda$ {\it be a left
Artinian ring and} $M\in \mod \Lambda$. {\it Then the following
statements are equivalent.}

(1) $M$ {\it is faithful.}

(2) $M$ {\it cogenerates every projective module.}

(3) $M$ {\it generates every injective module.}

\vspace{0.2cm}

{\bf Definition 2.3} ([AuB] or [EJ]) Let $\Lambda$ be a left and
right Noetherian ring. A module $M$ in $\mod \Lambda$ is called {\it
Gorenstein dimension zero} (or {\it Gorenstein projective}) if the
following conditions are satisfied: (1) $M$ is reflexive; (2) ${\rm
Ext}_{\Lambda}^i(M,\Lambda)=0={\rm
Ext}_{\Lambda^{op}}^i(M^*,\Lambda^{op})$ for any $i \geq 1$.

\vspace{0.2cm}

Recall that a module in $\mod \Lambda$ is called {\it
selforthogonal} if ${\rm Ext}^i_{\Lambda}(M, M)=0$ for any $i\geq
1$. Then it is trivial that {\bf GPC} is a special case of {\bf
GNC}.

\vspace{0.5cm}

\centerline{\large \bf 3. The Case for Artinian Algebras}

\vspace{0.2cm}

In this section, $\Lambda$ is an Artinian algebra. The following
lemma plays a crucial role in this section.

\vspace{0.2cm}

{\bf Lemma 3.1} {\it Let} $M\in \mod \Lambda$ {\it be an
indecomposable module. If there exists an exact sequence}
$M^t\rightarrow N\rightarrow 0$ {\it and} $\underline{{\rm
Hom}}_\Lambda (M^t, N)=0$ {\it for some} $t\geq 1$ {\it and} $N\in
\mod \Lambda$, {\it then} $M$ {\it is projective}.

\vspace{0.2cm}

{\bf Proof.} Let $(P(N),g)$ be the projective cover of $N$. Because
$\underline{{\rm Hom}}_\Lambda (M^t, N)=0$, we get a homomorphism
$h: M^t\rightarrow P(N)$ and the following commutative diagram with
exact rows:
$$\xymatrix{0\ar[r]&\Ker f\ar[r]\ar@{-->}[d]^{h'}&M^t\ar[r]
\ar[d]^{h}&N\ar[r]\ar@{=}[d]&0\\
0\ar[r]&\Ker g \ar[r]&P(N)\ar[r]^g&N\ar[r]&0 }
$$
where $h'$ is an induced homomorphism. Since $g$ is a superfluous
epimorphism, $h$ is epimorphic and splitable. So $P(N)$ is
isomorphic to a direct summand of $M^t$. Since $M$ is
indecomposable, $P(N)\cong M^s$ for some $s \geq 1$ and $M$ is
projective. \hfill{$\square$}

\vspace{0.2cm}

{\bf Lemma 3.2} {\it Let} $\Lambda$ {\it be a radical square zero
algebra and} $M\in \mod \Lambda$ {\it an indecomposable module. If}
$M$ {\it is torsionless and not simple, then} $M$ {\it is
projective}.

\vspace{0.2cm}

{\bf Proof.} Suppose $M\neq 0$. Then $M\neq J(\Lambda)M$ and there
exists a simple $\Lambda$-module $S$ such that
$M/J(\Lambda)M\rightarrow S\rightarrow 0$ is exact. Since $M$ is
indecomposable, we have a non-spilt epimorphism $f:M\rightarrow S$.

We claim that $\overline{{\rm Hom}}_\Lambda(M, S)=0$. If $S$ is
injective, then it is clear that $\overline{{\rm Hom}}_\Lambda(M,
S)=0$. If $S$ is not injective, then, by [AuR2, Proposition 4.3], we
have an almost spilt sequence $0\rightarrow S\rightarrow
E\rightarrow \Tr \mathbb{D}S\rightarrow 0$. Notice that
$J(\Lambda)^2=0$ by assumption, so $E$ is projective by [AuR2,
Proposition 5.7]. Since $M$ is not simple, ${\rm Ext}_\Lambda^1(\Tr
\mathbb{D}S,M)=0$ by [AuR2, Theorem 5.5]. So $\overline{{\rm
Hom}}_\Lambda(M, S)=0$ by Lemma 2.1. The claim is proved.

Since $M$ is torsionless,  there exists a projective $P\in\mod
\Lambda$ such that $0\rightarrow M\rightarrow P$ is exact. Then it
is easy to see that $\underline{{\rm Hom}}_\Lambda(M, S)=0$ and
there exists an exact sequence $M\rightarrow S\rightarrow 0$. By
Lemma 3.1, $M$ is projective. \hfill{$\square$}

\vspace{0.2cm}

{\bf Lemma 3.3} {\it Let} $\Lambda$ {\it be a local algebra with
radical square zero and} $M\in \mod \Lambda$ {\it an indecomposable
module. If} $M$ {\it is torsionless and} ${\rm
Ext^1_\Lambda}(M,M)=0$, {\it then} $M$ {\it is projective}.

\vspace{0.2cm}

{\bf Proof.} If $M$ is not simple, $M$ is projective by Lemma 3.2.
If $M$ is simple, then the condition ${\rm Ext^1_\Lambda}(M,M)=0$
implies $M$ is projective by [XC, Lemma 3]. \hfill{$\square$}

\vspace{0.2cm}

The following is the main result in this section.

\vspace{0.2cm}

{\bf Theorem 3.4} {\it Let} $\Lambda$ {\it be a local algebra with
radical square zero. Then a torsionless module} $M\in \mod \Lambda$
{\it is projective if} ${\rm Ext^1_\Lambda}(M,M)=0$.

\vspace{0.2cm}

{\bf Proof.} If $M\in \mod \Lambda$ is torsionless and ${\rm
Ext^1_\Lambda}(M,M)=0$, then $N$ is torsionless and ${\rm
Ext^1_\Lambda}(N,N)=0$ for any direct summand $N$ of $M$. Thus the
assertion follows immediately from Lemma 3.3. \hfill{$\square$}

\vspace{0.2cm}







In the following, we give some criteria for judging an
indecomposable torsionless module being projective.

\vspace{0.2cm}

{\bf Proposition 3.5} {\it Let} $M\in \mod \Lambda$ {\it be faithful
and indecomposable. Then} $M$ {\it is projective if} $M$ {\it is
torsionless}.

\vspace{0.2cm}

{\bf Proof.} By Lemma 2.1, for any $n\geq 1$, we have an
isomorphism:
$$\overline{{\rm Hom}}_{\Lambda^{op}} (\Lambda^{op}, \mathbb{D}(M^n))\cong {\rm Hom}
_{\underline{{\rm End}}(\Tr(M^n))} ({\rm Ext}^1_{\Lambda^{op}} (\Tr
(M^n), \Lambda^{op}), {\rm Ext}^1_{\Lambda^{op}} (\Tr (M^n),
\mathbb{D}(M^n)).$$ Notice that $M$ is torsionless, so ${\rm
Ext}^1_{\Lambda^{op}}(\Tr (M^n), \Lambda ^{op})\cong  {\rm
Ext}^1_{\Lambda^{op}}((\Tr M)^n, \Lambda ^{op})\cong ({\rm
Ext}^1_{\Lambda^{op}}(\Tr M,$ \linebreak $\Lambda ^{op})^n=0$, and
hence $\overline{{\rm Hom}}_{\Lambda^{op}} (\Lambda^{op},
\mathbb{D}(M^n))=0$ and $\underline{{\rm Hom}}_\Lambda (M^n,
\mathbb{D}\Lambda^{op})=0$. On the other hand, because $M$ is
faithful, by Lemma 2.2, there exists an $n \geq 1$ such that
$M^n\rightarrow \mathbb{D}\Lambda ^{op}\rightarrow 0$ is exact. So,
by Lemma 3.1, we have that $M$ is projective. \hfill{$\square$}

\vspace{0.2cm}

{\bf Proposition 3.6} {\it Let} $M\in \mod \Lambda$ {\it be faithful
and indecomposable. Then} $M$ {\it is a projective if} $\Tr M$ {\it
is torsionless (equivalently,} ${\rm Ext^1_{\Lambda}}(M,
\Lambda)=0${\it ) and} ${\rm Ext}^1_{\Lambda^{op}}(\Tr M,\Tr M)=0$.

\vspace{0.2cm}

{\bf Proof.} Since $\Tr M$ is torsionless, there exists a
monomorphism $0\rightarrow \Tr M\rightarrow (\Lambda ^{op})^n$. Then
$\mathbb{D}(\Lambda ^{op})^n\rightarrow \mathbb{D}\Tr M\rightarrow
0$ is exact. Because $M$ is faithful, there exists an $m\geq 1$ such
that $M^m\rightarrow \mathbb{D}({\Lambda^{op}})^n\rightarrow 0$ is
exact. So we have an exact sequence $M^m\rightarrow \mathbb{D}\Tr
M\rightarrow 0$. On the other hand, Since ${\rm
Ext}^1_{\Lambda^{op}}(\Tr (M^m),\Tr M)\cong {\rm
Ext}^1_{\Lambda^{op}}((\Tr M)^m,\Tr M)\cong ({\rm
Ext}^1_{\Lambda^{op}}(\Tr M,\Tr M))^m$ \linebreak $=0$ by
assumption, $\overline{{\rm Hom}}_{\Lambda^{{op}}}(\Tr
M,\mathbb{D}(M^m))=0$ by Lemma 2.1. So $\underline{{\rm
Hom}}_{\Lambda}(M^m,\mathbb{D}\Tr M)=0$, and hence $M$ is projective
by Lemma 3.1. \hfill{$\square$}

\vspace{0.2cm}

Recall from [AuR1] that the generalized Nakayama conjecture ({\bf
GNC}) states that a module $M\in \mod \Lambda$ is projective if
${\rm Ext}_\Lambda ^i(M\bigoplus \Lambda, M\bigoplus \Lambda)=0$ for
any $i \geq 1$, which still remains open. The following result
provides some support to this conjecture.

\vspace{0.2cm}

{\bf Proposition 3.7} {\it Let} $S$ {\it be a faithful and simple
in} $\mod \Lambda$. {\it If} ${\rm Ext_\Lambda ^1}(S\bigoplus
\Lambda,S\bigoplus \Lambda)=0$, {\it then} $S$ {\it is projective}.

\vspace{0.2cm}

{\bf Proof.} Since ${\rm Ext^1_{\Lambda}}(S, S)=0$ by assumption,
$\overline{{\rm Hom}}_\Lambda(S,\mathbb{D}\Tr S)=0$ by Lemma 2.1. So
$\underline{{\rm Hom}}_{\Lambda^{op}}(\Tr S, \mathbb{D}S)=0$.

If ${\rm Hom}_\Lambda(S,\mathbb{D}\Tr S)\neq 0$, then we have an
epimorphism $\Tr S\rightarrow \mathbb{D}S\rightarrow 0$ since $S$ is
simple. By Lemma 3.1, $\Tr S$ is projective and $\Tr S=0$. So $S$ is
projective.

If ${\rm Hom}_\Lambda(S,\mathbb{D}\Tr S)=0$, then ${\rm
Hom}_\Lambda(S^m,\mathbb{D}\Tr S)=0$ for any $m \geq 1$. Because
${\rm Ext^1_{\Lambda}}(S,\Lambda )=0$ by assumption, $\Tr S$ is
torsionless and there exists a monomorphism $0\rightarrow \Tr
S\rightarrow (\Lambda^{op})^n$. So
$\mathbb{D}({\Lambda^{op}})^n\rightarrow \mathbb{D}\Tr S\rightarrow
0$ is exact. Because $S$ is faithful, there exists an $m \geq 1$
such that $S^m\rightarrow \mathbb{D}({\Lambda^{op}})^n\rightarrow 0$
is exact. So we have an epimorphism $S^m\rightarrow \mathbb{D}\Tr
S\rightarrow 0$. It implies that $\mathbb{D}\Tr S=0$ and $\Tr S=0$.
Thus $S$ is projective. \hfill{$\square$}

\vspace{0.5cm}

\centerline{\large \bf 4. The Case for Commutative Artinian Rings}

\vspace{0.2cm}

In this section, $\Lambda$ is a commutative Artinian ring. According
to the localization theory of commutative ring, by Theorem 3.4, we
have the following

\vspace{0.2cm}

{\bf Theorem 4.1} {\it If} $\Lambda$ {\it is radical square zero,
then a torsionless module} $M\in \mod \Lambda$ {\it is projective
if} ${\rm Ext^1_\Lambda}(M,M)=0$.

\vspace{0.2cm}

Let $M$ and $N$ be in $\mod \Lambda$. We define a homomorphism
$\zeta : M\otimes _{\Lambda}N\to {\rm Hom}_\Lambda(M^*, N)$ of
$\Lambda$-modules by $\zeta (m\otimes n)(g)=g(m)n$ for any $m\otimes
n\in M\otimes_\Lambda N$ and $g\in M^*$. Then we obtain a natural
transformation $\zeta (-):M\otimes_\Lambda -\to {\rm Hom} _\Lambda
(M^*,-)$ of functors from $\mod \Lambda$ to itself.

\vspace{0.2cm}

{\bf Lemma 4.2} ([AuB, Proposition 2.6]) {\it For any} $M\in\mod
\Lambda$, {\it there exists an exact sequence of functors from}
$\mod \Lambda$ {\it to itself:}
$$0 \to {\rm Ext_\Lambda ^1(Tr}M,-) \to M\otimes _{\Lambda}-
\buildrel {{\zeta (-)}} \over \longrightarrow {\rm
Hom_\Lambda}(M^*,-) \to  {\rm Ext_\Lambda ^2(Tr}M,-) \to 0.$$

\vspace{0.2cm}

{\bf Definition 4.3} ([AuR3]) Assume that $\mathscr{X}$ is a full
subcategory of $\mod \Lambda$ and $Y\in \mod \Lambda$,
$X\in\mathscr{X}$. The morphism $f:X\to Y$ is said to be a {\it
right} $\mathscr{X}$-{\it approximation} of $Y$ if ${\rm
Hom}_{\Lambda}(X', X)\to {\rm Hom}_{\Lambda}(X', Y)\to 0$ is exact
for any $X'\in\mathscr{X}$. The morphism $f:X\to Y$ is said to be
{\it right minimal} if an endomorphism $g:X\to X$ is an automorphism
whenever $f=fg$. The subcategory $\mathscr{X}$ is said to be {\it
contravariantly finite} in $\mod \Lambda$ if every $Y\in\mod
\Lambda$ has a right $\mathscr{X}$-approximation. The notions of
{\it (minimal) left} $\mathscr{X}$-{\it approximations} and {\it
covariantly finite subcategories} of $\mod \Lambda$ may be defined
dually. The subcategory $\mathscr{X}$ is said to be {\it
functorially finite} in $\mod \Lambda$ if it is both contravariantly
finite and covariantly finite in $\mod \Lambda$.

\vspace{0.2cm}

For a module $M \in \mod \Lambda$, we denote ${^{\perp_1}M} =\{X\in
\mod \Lambda |{\rm Ext}_\Lambda ^1(X, M)=0\}$.

\vspace{0.2cm}

{\bf Lemma 4.4} ([T, Lemma 6.9]) {\it Let} $M\in \mod \Lambda$ {\it
with} $M\in {^{\perp_1}M}$. {\it Then for any} $N\in \mod \Lambda$,
{\it there exists an exact sequence} $0\to F\to E\to N\to 0$ {\it
with} $F=M^{(n)}$ {\it and} $E\in {^{\perp_1}M}$, {\it where} $n$
{\it is number of generators of} ${\rm Ext}_\Lambda^1(N,M)$ {\it as
an} ${\rm End}(M)$-{\it module. Hence} ${^{\perp_1}M}$ {\it is
contravariantly finite.}

\vspace{0.2cm}

{\bf Lemma 4.5} ([AuS, Proposition 7.1]) $^{\bot _1}\Lambda$ {\it is
functorially finite in} $\mod \Lambda$.

\vspace{0.2cm}

Now we give the main result in this section.

\vspace{0.2cm}

{\bf Theorem 4.6} {\it A torsionless module} $M\in \mod \Lambda$
{\it is projective if the following conditions are satisfied}:

(1) ${\rm Ext}^i_{\Lambda}(M,\Lambda)=0$ {\it for} $i=1,2,3$.

(2) ${\rm Ext}^i_{\Lambda}(M,M)=0$ {\it for} $i=1,2$.

\vspace{0.2cm}

{\bf Proof.} Without loss of generality, we can assume that
$\Lambda$ is local with unique maximal ideal $m$ and residue field
$k(=R/m)$.

From Lemma 4.4, we know there exists a right $^{\perp_1}
M$-approximation: $0\rightarrow M^n\rightarrow E^\prime\rightarrow
k\rightarrow 0$ for the simple $\Lambda$-module $k$, where $n$ is
the number of the generators of ${\rm Ext}_\Lambda ^1(k,M)$ as an
${\rm End}(M)$-module. If $n=0$, then ${\rm Ext}_\Lambda ^1(k,M)=0$.
So $M$ is injective. But $M$ is torsionless by assumption, thus $M$
is projecitve.

Now suppose $n\geq 1$. Consider the minimal right $^{\perp_1}
M$-approximation of $k$: $0\rightarrow M^m\rightarrow E\rightarrow
k\rightarrow 0$. By applying the functor $\Tr M\otimes_\Lambda -$ to
it, we obtain a commutative diagram with exact rows:
$$\xymatrix{
&\Tr M\otimes_\Lambda M^m\ar[r]\ar[d]^{\zeta (M^m)}& \Tr
M\otimes_\Lambda E\ar[r]\ar[d]^{\zeta(\Tr M)}&\Tr M\otimes_\Lambda
k\ar[d]^{\zeta(k)}\ar[r]&0\\
0\ar[r]&{\rm Hom}_\Lambda ((\Tr M)^*, M^m)\ar[r]^\alpha&{\rm
Hom}_\Lambda ((\Tr M)^*, E)\ar[r]^\beta&{\rm Hom}_\Lambda ((\Tr
M)^*, k)& }
$$
Since ${\rm Ext}_\Lambda ^i(M,M)=0$ for $i=1,2$, ${\zeta(M^m)}$ is
an isomorphism by Lemma 4.2.

Consider the homomorphism $\zeta(k):\Tr M\otimes_{\Lambda} k\to{\rm
Hom}_\Lambda ((\Tr M)^*, k)$ via $\zeta(k)(a\otimes
\overline{r})(f)=f(a)\overline{r}$ for any $a\in \Tr M, f\in (\Tr
M)^*$ and $\overline{r}\in k$. Because $\Tr M$ has no projective
summands, $f$ is not epimorphic. Notice that $(\Lambda,m,k)$ is
local, $f(\Tr M)\subseteq m$. It follows that $\zeta(k)(a\otimes
\overline{r})(f)=f(a)\overline{r}=0$ and $\zeta(k)=0$. Then we have
that $\beta\circ{\zeta(\Tr M)}=0$, and thus there exists a
homomorphism $\gamma :\Tr M\otimes _{\Lambda}E\rightarrow {\rm
Hom}_\Lambda ((\Tr M)^*, M^m)$ such that $\alpha \circ \gamma =
{\zeta(\Tr M)}$. Since $\alpha$ is monomorphic, the sequence
$0\to\Tr M\otimes _{\Lambda} M^m\to \Tr M\otimes _{\Lambda}E\to \Tr
M\otimes _{\Lambda} k\to 0$ (the upper row in the above diagram) is
exact and split. Then we get a commutative diagram with exact rows:
$$\xymatrix{
0\ar[r]&(\Tr M\otimes _{\Lambda}k)^*\ar[r]\ar[d]^{\cong}& (\Tr
M\otimes  _{\Lambda}E)^*\ar[r]\ar[d]^{\cong}&(\Tr M\otimes
_{\Lambda}
M^m)^*\ar[d]^{\cong}\ar[r]&0\\
0\ar[r]&{\rm Hom}_\Lambda (k, (\Tr M)^*)\ar[r]&{\rm Hom}_\Lambda (E,
(\Tr M)^*)\ar[r]&{\rm Hom}_\Lambda (M^m,(\Tr M)^*)\ar[r]&0 }
$$
By the exactness of the bottom row in the above diagram, we have an
exact sequence $0\to {\rm Ext}_\Lambda^1(k, (\Tr M)^*)\to{\rm
Ext}_\Lambda^1(E, (\Tr M)^*)$. By the claim below, ${\rm
Ext}_\Lambda^1(E, (\Tr M)^*)=0$, so ${\rm Ext}_\Lambda^1(k, (\Tr
M)^*)=0$ and $(\Tr M)^*$ is injective. Notice that $\Ker f\cong (\Tr
M)^*$ in the minimal projective resolution $P_1 \buildrel {f} \over
\to P_0 \to M \to 0$, so the projective dimension of $M$ is at most
1. On the other hand, ${\rm Ext}^1_{\Lambda}(M, \Lambda)=0$ by
assumption, then it is easy to see that $M$ is projective.

{\bf Claim}: ${\rm Ext}_\Lambda^1(E, (\Tr M)^*)=0$.

Consider the exact sequence $0\to (\Tr M)^* \to P_1\to P_0\to M\to
0$. Since ${\rm Ext}^i_{\Lambda}(M,\Lambda)=0$ for $i=1,2,3$, $(\Tr
M)^*\in {^{\perp_1}\Lambda}$ . Let $0\to (\Tr M)^* \to
Z\stackrel{h}{\to} E\to 0$ be any exact sequence in ${\rm
Ext}_\Lambda ^1(E,(\Tr M)^*)$. Consider the following pullback
diagram:
$$\xymatrix{
&&0\ar[d]&0\ar[d]&\\
0\ar[r]&(\Tr M)^*\ar[r]\ar@{=}[d]&X\ar[d]\ar[r]&M^m\ar[r]\ar[d]&0\\
0\ar[r]&(\Tr M)^*\ar[r]&Z\ar[d]\ar[r]&E\ar[r]\ar[d]&0\\
&&k\ar[d]\ar@{=}[r]&k\ar[d]&\\
&&0&0&}
$$
Since $^{\perp_1} \Lambda$ is closed under extensions, $X$ is in
$^{\perp_1}\Lambda$. On the other hand, $^{\perp_1}\Lambda$
covariantly finite by Lemma 4.5, ${\rm
Hom}_{\Lambda}(M^m,-)|_{^{\perp_1}\Lambda}\to {\rm
Ext}_\Lambda^1(k,-)|_{^{\perp_1}\Lambda}\to 0$ is exact. So there
exists a commutative diagram with exact rows:
$$\xymatrix{
0\ar[r]&M^m\ar[d]\ar[r]&E\ar[r]\ar[d]&k\ar@{=}[d]\ar[r]&0\\
0\ar[r]&X\ar[r]&Z\ar[r]&k\ar[r]&0 }
$$
Then we obtain a commutative diagram with exact rows:
$$\xymatrix{
0\ar[r]&M^m\ar[d]\ar[r]&E\ar[r]\ar[d]&k\ar@{=}[d]\ar[r]&0\\
0\ar[r]&X\ar[r]\ar[d]&Z\ar[r]\ar[d]^h&k\ar@{=}[d]\ar[r]&0\\
0\ar[r]&M^m\ar[r]&E\ar[r]&k\ar[r]&0 }
$$
Because $0\rightarrow M^m\rightarrow E\rightarrow k\rightarrow 0$ is
right minimal, it follows that the composition $M^m\to X\to M^m$ is
an isomorphism. Thus the composition $E\to Z\stackrel{h}{\to} E$ is
also an isomorphism, and therefore the exact sequence $0\to (\Tr
M)^* \to Z\stackrel{h}{\to} E\to 0$ splits and ${\rm
Ext}_\Lambda^1(E, (\Tr M)^*)=0$.  \hfill{$\square$}

\vspace{0.2cm}

As an immediate consequence of Theorem 4.6, we get the following
result, which means that {\bf GPC} is true for commutative Artinian
rings. This also provides some support to {\bf GNC}.

\vspace{0.2cm}

{\bf Theorem 4.7} {\it A Gorenstein projective module in} $\mod
\Lambda$ {\it is projective if and only if it is selforthogonal}.

\vspace{0.5cm}

{\bf Acknowledgements} This research was partially supported by the
Specialized Research Fund for the Doctoral Program of Higher
Education (Grant No. 20060284002), NSFC (Grant No. 10771095) and NSF
of Jiangsu Province of China (Grant No. BK2007517).

\end{document}